\documentclass[a4paper,reqno]{amsart}
\usepackage{amsmath,amsfonts,amscd,amssymb,amsthm,graphicx,ifthen}
\newcommand{\nz}{\mathbb{N}}       
\newcommand{\rz}{\mathbb{R}}       
	\theoremstyle{plain}
	\newtheorem{SATZ}{Proposition}[section]   
	\newtheorem{THM}[SATZ]{Theorem}
	\newtheorem{LEM}[SATZ]{Lemma}

	\theoremstyle{definition}

	\newtheorem{BEM}[SATZ]{Remark}
	  \newenvironment{BEW}[1][Proof]{\begin{proof}[#1]}{\end{proof}}
	  
\numberwithin{figure}{section}
\numberwithin{table}{section}
\numberwithin{equation}{section}
\newcommand{\dx}{{\rm d}}
\DeclareMathOperator{\Span}{span}
\DeclareMathOperator{\Div}{div}
\DeclareMathOperator{\diam}{diam}
\DeclareMathOperator*{\pinfp}{\phantom{p}inf\phantom{p}}
\DeclareMathOperator*{\esssup}{ess.sup}
\newcommand{\hdiamE}[1][E]{h_{#1}}
\newcommand{\hdiamK}[1][K]{h_{#1}}

\newcommand{\mdiam}[1]{\hslash_{#1}}
\DeclareMathOperator{\shp}{C}

\newcommand{\shape}[1][\CT]{\shp_{#1}}
\newcommand{\halb}{\frac 12}
\newcommand{\cd}{\Va \cdot \nabla}
\newcommand{\cdm}{\bar{\Va}_{\CM}\cdot\nabla}
\newcommand{\cdt}{\bar{\Va}_{\CT}\cdot\nabla}
\newcommand{\al}{\alpha}
\newcommand{\be}{\beta}

\newcommand{\Ga}{\Gamma}

\newcommand{\De}{\Delta}
\newcommand{\ep}{\varepsilon}
\newcommand{\ka}{\kappa}

\newcommand{\om}{\omega}
\newcommand{\Om}{\Omega}
\newcommand{\vfi}{\varphi}
\newcommand{\si}{\sigma}
\newcommand{\Si}{\Sigma}

\newcommand{\rN}{\rvert}

\newcommand{\betrag}[2][]{\left\lvert #2 \right\rvert_{#1}} 
\newcommand{\norm}[2][]{\left\lVert #2 \right\rVert_{#1}} 
\newcommand{\normE}[2][]{\left\lVert{\hskip -0.25em} \left\lvert #2 \right\rVert {\hskip -0.25em} \right\rvert_{#1}} 
\newcommand{\normD}[2][]{\left\lVert{\hskip -0.25em} \left\lvert #2 \right\rVert {\hskip -0.25em} \right\rvert_{\ast#1}} 
\newcommand{\normL}[3][2]{				 
	{
			\ifthenelse{\equal{#1}{2}}
				{
					\ifthenelse{\equal{#2}{\Om} \or \equal{#2}{\Omega}}
					{
						\left\lVert #3 \right\rVert
					}
					{
						\left\lVert #3 \right\rVert_{#2}
					}
				}
				{
					\ifthenelse{\equal{#2}{\Om} \or \equal{#2}{\Omega}}
					{
						\left\lVert #3 \right\rVert_{#1}
					}
					{
						\left\lVert #3 \right\rVert_{#1;#2}
					}
				}
		}
	}
\newcommand{\normH}[3][1]{				 
	{
			\ifthenelse{\equal{#2}{\Om} \or \equal{#2}{\Omega}}
			{
				\left\lVert #3 \right\rVert_{#1}
			}
			{
				\left\lVert #3 \right\rVert_{#1;#2}
			}
		}
	}
\newcommand{\normW}[4][1]{				 
	{
			\ifthenelse{\equal{#3}{\Om} \or \equal{#3}{\Omega}}
			{
				\left\lVert #4 \right\rVert_{#1,#2}
			}
			{
				\left\lVert #4 \right\rVert_{#1,#2;#3}
			}
		}
	}
\newcommand{\paar}[3][]{\left\langle #2\,,\,#3 \right\rangle_{#1}} 
\newcommand{\jump}[2][E]{{\mathbb J}_{#1}(#2)}
\newcommand\Va{\mathbf{a}}

\newcommand\Vn{\mathbf{n}}

\newcommand\CE{\mathcal{E}}

\newcommand\CI{\mathcal{I}}

\newcommand\CM{\mathcal{M}}

\newcommand\CT{\mathcal{T}}

\title[Stabilized Finite Element Methods]{Robust A Posteriori Error Estimates for Stabilized Finite Element Methods}
\author{L. Tobiska}
\author{R. Verf{\"u}rth}
\address{Otto von Guericke Universit{\"a}t Magdeburg \\Fakult{\"a}t f{\"u}r Mathematik \\Institut f{\"u}r Analysis und Numerik \\ Universit{\"a}tsplatz 2 \\D-39106 Magdeburg \\Germany}
\address{Ruhr-Uni\-ver\-si\-t{\"a}t Bo\-chum \\ Fa\-kul\-t{\"a}t f{\"u}r Ma\-the\-ma\-tik \\ D-44780 Bo\-chum \\ Germany}
\email{tobiska@ovgu.de}
\email{ruediger.verfuerth@rub.de}
\date{\today}
\keywords{robust a posteriori error estimates, stabilized finite element methods, streamline diffusion methods, local projection schemes, continuous interior penalty methods, subgrid-scale methods, stationary convection-diffusion equations, non-stationary convection-diffusion equations}
\subjclass{65N30, 65N15, 65J10}
\date{\today}
\begin{document}
\begin{abstract}
There is a wide range of stabilized finite element methods for stationary and non-stationary convection-diffusion equations such as streamline diffusion methods, local projection schemes, subgrid-scale techniques, and continuous interior penalty methods to name only a few. We show that all these schemes give rise to the same robust a posteriori error estimates, i.e.  the multiplicative constants in the upper and lower bounds for the error are independent of the size of the convection or reaction relative to the diffusion. Thus, the same error indicator can be used modulo higher order terms caused by data approximation.  
\end{abstract}
\maketitle
%
%
\section{Introduction}\label{P:intro}
There is a wide range of stabilized finite element methods for stationary and non-stationary convection-diffusion problems such as streamline diffusion methods (cf.\ eg.\ \cite{FraFreHug,FraJohMatTob,HB79,Nae82,RST08}), local projection schemes (cf.\ eg.\ \cite{AraPozVal13,BB01,BB04,BB06,HT12,KT11,KT13,MST07,MT13,TW10}), subgrid-scale techniques (cf.\ eg.\ \cite{AchBenFatSouWar12,AchFatSou09,ErnGue04,Gal11,Gue99,Gue01,GMQ06,ZhaHe10}), and continuous interior penalty methods (cf.\ eg.\ \cite{AEB07,BraBurJohLub,Bur05,Bur07,BE07,BH04,DEV13,EG13,YanZho09}) to name only a few. In this article we show that all these schemes give rise to the same robust a posteriori error estimates. Here, as usual, robustness means that the upper and lower bounds for the error are uniform with respect to the size of convection or reaction terms relative  to the diffusion. Our analysis is based on the general approach of \cite{Ver13} which gives the generic robust equivalence of error and residual an
 d provides robust global upper and local lower bounds for the residual up to a consistency error in the upper bound. The latter depends on the particular stabilization method. Consequently, our main task consists in deriving explicit and computable upper bounds for the consistency errors of the various schemes. This is the subject of Lemmas \ref{L:sdfem} -- \ref{L:cip} below, where the result for the streamline diffusion method, Lemma \ref{L:sdfem}, is a reformulation of known results (cf.\ \cite[\S\S 4.4, 6.2]{Ver13} and the references given there). The main results of this article are Theorems \ref{T:error-est} and \ref{T:nonstat-error-est} below which provide robust residual a posteriori error estimates for stationary and non-stationary convection-diffusion equations, respectively.

There are also other proposals in the literature for estimating the error with respect to an a-priori given mesh dependent norm 
(cf.\ eg.\ \cite{AABR13,JN13}) or even control certain functionals of the solution like in the dual weighted residual approach 
(cf.\ eg.\ \cite{BR03}). Note, however, that these methods often use assumptions which are difficult to establish in practice. 
For a posteriori error control for other classes of problems we finally refer to \cite{AO00,CEHL12,Ver13}.
%
%
\section{Stationary Convection-Diffusion Equations}\label{P:stat-problem}
\subsection{Variational Problem}\label{P:stat-var-pb}
In this section, we consider the stationary convection-diffusion equation
\begin{equation}\label{E:stat-conv-diff}
\begin{aligned}
-\ep \Delta u + \cd u + b u &=  f &&\text{in } \Om \\
u &= 0 &&\text{on } \Ga_D \\
\ep \frac{\partial u}{\partial n} &= g &&\text{on } \Ga_N
\end{aligned}
\end{equation}
in a polygonal domain $\Om$ in $\rz^d$, $d \ge 2$, with Lipschitz boundary $\Ga$ consisting of two disjoint components $\Ga_D$ and $\Ga_N$. We assume that the data satisfy the following conditions:
\begin{enumerate}
\renewcommand{\labelenumi}{(S\arabic{enumi})}
\item $0 < \ep \ll 1$,
\item $\Va \in W^{1,\infty}(\Om)^d$, $b \in L^\infty(\Om)$,
\item there are two constants $\be \ge 0$ and $c_b \ge 0$, which do not depend on $\ep$, such that $-\halb \Div \Va + b \ge \be$ and $\normL[\infty]{\Om}{b} \le c_b \be$,
\item the Dirichlet boundary $\Ga_D$ has positive
$(d-1)$-dimensional Hausdorff measure and includes the inflow boundary, i.e.\ $\left\{x \in \Ga: \Va(x) \cdot \Vn(x) < 0\right\} \subset \Ga_D$.
\end{enumerate}
Assumption (S3) allows us to handle simultaneously the case of a non-vanishing reaction term and the one of absent reaction. If $b= 0$ we set $\be=c_b = 0$. Assumption (S1) of course means that we are interested in the convection-dominated regime.

For the variational formulation of problem \eqref{E:stat-conv-diff} we denote by $H^1_D(\Om)$ the standard Sobolev space of all functions in $L^2(\Om)$ having their weak first order derivatives in $L^2(\Om)$ and vanishing on $\Ga_D$ in the sense of traces and define a bilinear form $B$ on $H^1_D(\Om) \times H^1_D(\Om)$ and a linear functional $\ell$ on $H^1_D(\Om)$ by
\begin{equation}\label{E:def-B-ell}
B(u, v) = \int_\Om \left\{ \ep \nabla u \cdot \nabla v + \cd u v + b u v \right\}, \quad \paar{\ell}{v} = \int_\Om f v + \int_{\Ga_N} g v. 
\end{equation}
The variational problem then consists in finding $u \in H^1_D(\Om)$ such that
\begin{equation}\label{E:stat-var-pb}
B(u,v) = \paar{\ell}{v}
\end{equation}
holds for all $v \in H^1_D(\Om)$.

The well-posedness of problem \eqref{E:stat-var-pb} and the robustness of the a posteriori error estimates hinges on a proper choice of norms. More specifically, we denote by
\begin{equation*}
\normE{v}  = \left\{ \ep \normL{\Om}{\nabla v}^2 + \be \normL{\Om}{v}^2 \right\}^\halb
\end{equation*}
the energy norm associated with symmetric part of $B$ and by
\begin{equation*}
\normD{\vfi} = \sup_{v \in H^1_D(\Om) \setminus \{ 0 \}} \frac{\paar{\vfi}{v}}{\normE{v}},
\end{equation*}
the corresponding dual norm on $H^{-1}(\Om) = \left( H^1_D(\Om) \right)^\ast$. Here, $\normL{\om}{\cdot}$ is the standard $L^2$-norm on any measurable subset $\om$ of $\Om$ and $\normL{\Om}{\cdot} = \norm[\Om]{\cdot}$. With this choice of norms we have for all $v, w \in H^1_D(\Om)$ \cite[Prop. 4.17]{Ver13}
\begin{equation*}
\begin{split}
B(v,v) &\ge \normE{v}^2, \\
B(v,w) &\le \max\left\{c_b , 1\right\} \left\{\normE{v} + \normD{\cd v} \right\} \normE{w}
\end{split}
\end{equation*}
and
\begin{equation*}
\pinfp_{v \in H^1_D(\Om) \setminus \{0\}} \sup_{w \in H^1_D(\Om)
\setminus \{0\}} \frac{B(v ,w)}{\left\{\normE{v} + \normD{\cd v} \right\} \normE{w}} \ge \frac{1}{2 + \max\left\{c_b, 1\right\}}.
\end{equation*}
This in particular implies that problem \eqref{E:stat-var-pb} admits for every right-hand side $\ell \in H^{-1}(\Om)$ a unique solution $u \in H^1_D(\Om)$ and that
\begin{equation}\label{E:robustness}
c_\sharp \normD{\ell} \le \normE{u} + \normD{\cd u} \le c^\sharp \normD{\ell}
\end{equation}
with constants $c_\sharp$ and $c^\sharp$ only depending on $c_b$ and independent of $\ep$ and $\be$.

\subsection{Discretization}\label{P:stat-discr-pb}
For the discretization of problem \eqref{E:stat-conv-diff}, we denote by $\CT$ a partition of $\Om$ which satisfies the following conditions.
\begin{itemize}
\item The closure of $\Om$ is the union of all elements in $\CT$.
\item The Dirichlet boundary $\Ga_D$ is the union of $(d-1)$-dimensional faces of elements in $\CT$.
\item Every element has at least one vertex in $\Om \cup \Ga_N$.
\item Every element in $\CT$ is either a simplex or a parallelepiped, i.e.\ it is the image of the $d$-dimensional reference simplex $\widehat K_d = \left\{ x\in \rz^d : x_1 \ge 0, \, \ldots, \, x_d \ge 0, \right.$\linebreak[4]$\left.x_1 + \ldots + x_d \le 1 \right\}$ or of the $d$-dimensional reference cube $\widehat K_d = [0, 1]^d$ under an affine mapping (\emph{affine-equivalence}).
\item Any two elements in $\CT$ are either disjoint or share a complete lower dimensional face of their boundaries (\emph{admissibility}).
\item For any element $K$, the ratio of its diameter $\hdiamK$ to the diameter $\rho_K$ of the largest ball inscribed into $K$ is bounded independently of $K$ (\emph{shape-regularity}).
\end{itemize}
As a measure for the shape-regularity we set as usual
\begin{equation}\label{E:shape-param}
\shape = \max_{K \in \CT} \frac{\hdiamK}{\rho_K}.
\end{equation}
The set of all $(d-1)$-dimensional faces of elements in $\CT$ is denoted by $\CE$. An additional subscript $\Om$, $\Ga_N$, or $K$ to $\CE$ indicates that only those faces are considered that are contained in the corresponding set. The union of all faces is called the skeleton of $\CT$ and denoted by $\Si$.
\\
As usual, we associate with every face $E \in \CE$ a unit vector $\Vn_E$ which is orthogonal to $E$ and which points to the outside of $\Om$ if $E$ is a face on the boundary $\Ga$. Finally, $\jump{\cdot}$ denotes the jump across $E$ in direction $\Vn_E$. Note, that $\jump{\cdot}$ depends on the orientation of $\Vn_E$ but that expressions of the form $\jump{\Vn_E \cdot \nabla v}$ are independent thereof. 

For every multi-index $\al \in \nz^d$, we set for abbreviation
\begin{equation*}
\betrag[1]{\al} = \al_1 + \ldots + \al_d, \; \betrag[\infty]{\al} = \max \left\{\al_i : 1 \le i \le d \right\}, \; x^\al = x_1^{\al_1} \cdot \ldots \cdot x_d^{\al_d}.
\end{equation*}
With every integer $k$ we then associate the standard spaces $P_k(\widehat K_d)$ and $Q_k(\widehat K_d)$ of polynomials by
\begin{equation*}
\begin{aligned}
P_k(\widehat K_d) &=
\Span \left\{x^\al : \betrag[1]{\al} \le k \right\} &&\text{for the reference simplex}, \\
Q_k(\widehat K_d)&=\Span \left\{x^\al : \betrag[\infty]{\al} \le k\right\} &&\text{for the reference cube}
\end{aligned}
\end{equation*}
and set for every element $K$ in $\CT$ and for $S\in\{P,Q\}$
\begin{equation*}
S_k(K) = \left\{ \vfi \circ F_K^{-1} : \vfi \in S_k(\widehat K_d) \right\},
\end{equation*}
where $F_K$ is an affine diffeomorphism from $\widehat K_d$ onto $K$. Using this notation, we define finite element spaces by
\begin{equation*}
\begin{split}
S^{k,-1}(\CT) &= \left\{ \vfi : \Om \to \rz : \vfi \rN_K \in S_k(K) \;\text{ for all } K \in \CT \right\}, \\
S^{k,0}(\CT) &= S^{k,-1}(\CT) \cap C(\overline \Om), \\
S_D^{k,0}(\CT) &= S^{k,0}(\CT) \cap H_D^1(\Om) = \left\{ \vfi \in S^{k,0}(\CT) : \vfi = 0 \,\,\text{ on } \Ga_D \right\}.
\end{split}
\end{equation*}
The number $k$ may be $0$ for the first space, but must be at least $1$ for the other spaces. Notice in particular that $P^{k,-1}(\CT)$ and $P^{k,0}(\CT)$ consist of piecewise polynomials of total degree at most $k$, that $Q^{k,-1}(\CT)$ and $Q^{k,0}(\CT)$ consist of piecewise polynomials of maximal degree at most $k$, and that $P^{k,-1}(\CT) \subset Q^{k,-1}(\CT)$ and $P^{k,0}(\CT) \subset Q^{k,0}(\CT)$.

Beside finite element spaces on the partition $\CT$ we consider spaces $S^{k,-1}(\CM)$, $S\in\{P,Q\}$ associated with a macro-partition $\CM$ subordinate to $\CT$. Note that for mapped reference cubes $M\in\CM$ we also consider discontinuous spaces $P^{k,-1}(\CM)$.  Moreover, we in addition assume that the partition $\CM$ either equals $\CT$ or is a coarsening of $\CT$ such that the elements of $\CT$ and $\CM$ are of comparable size, i.e. $\max_{M \in \CM} \max_{K \in \CT; K \subset M} \frac{\hdiamK[M]}{\hdiamK}$ is of moderate size independently of $\ep$ and $\be$.

Our discrete solution space is $V(\CT)$ with $S^{k,0}_D(\CT)\subset V(\CT)\subset S^{l,0}_D(\CT)$, 
$l\ge k$. The discretization of problem \eqref{E:stat-conv-diff} then consists in finding 
$u_\CT \in V(\CT)$ such that
\begin{equation}\label{E:discr-stat-pb}
\begin{split}
B(u_\CT,v_\CT) + S_\CT(u_\CT,v_\CT) = \paar{\ell}{v_\CT}
\end{split}
\end{equation}
holds for all $v_\CT \in V(\CT)$. Here, the term $S_\CT$ specifies the particular stabilization. It is supposed to be linear in its second argument and affine in its first argument. Note that $S_\CT$ may contain contributions of the data $f$ and $g$. Of course, the choice $S_\CT = 0$ is also possible and corresponds to a standard finite element method without stabilization. In the following subsections we consider in some more detail the stabilized schemes that are at the focus of this article. We always assume that the discrete problem \eqref{E:discr-stat-pb} admits a unique solution $u_\CT$. For the schemes below this is proved in the references given below by establishing the coercivity of the bilinear form $v,w \mapsto B(v,w) + S_\CT(v,w) - S_\CT(0,w)$ with respect to a suitable mesh-dependent norm.

\subsubsection{Streamline diffusion method}\label{SDFEM} This residual based stabilization 
method was introduced in \cite{HB79} and analyzed starting with \cite{Nae82}
under different aspects in a large number of articles, for an overview see eg. \cite{RST08}.
Only one partition $\CM=\CT$ of $\Omega$ is considered. The stabilizing term has the form
\begin{equation*}
S_\CT(u_\CT,v_\CT) = \sum_{K \in \CT} \vartheta_K \int_K \left\{-\ep \Delta u_\CT + \cd u_\CT + b u_\CT - f\right\} \cd v_\CT
\end{equation*}
with
\begin{equation}\label{E:cs-sdfem}
\vartheta_K \normL[\infty]{K}{\Va} \le c_S \hdiamK  \quad \text{for all }K \in \CT.
\end{equation}

\subsubsection{Local projection scheme}\label{LPS} This stabilization method has been first introduced for equal order interpolations of the Stokes problem in \cite{BB01}, extended to the transport problem in \cite{BB04}, and analyzed for the Oseen problem in \cite{BB06,MST07,MT13}. There are different versions on the market \cite{AraPozVal13,HT12,KT11,KT13,TW10}, here we consider the one-level approach ($\CT=\CM$) and the two-level approach ($\CT$ a subdivision of $\CM$) with two types of stabilizing terms, controlling the fluctuations of the derivatives in streamline direction
\begin{equation*}
S_\CT(u_\CT,v_\CT) = \sum_{M\in \CM} \vartheta_M \int_M \ka_\CM \left( \cdm u_\CT \right)  \ka_\CM \left( \cdm v_\CT \right)
\end{equation*}
with
\begin{equation}\label{E:cs-lps1}
\vartheta_M \normL[\infty]{M}{\Va}\le c_S \hdiamK[M]\quad \text{for all }M \in \CM
\end{equation}
or the fluctuations of the full gradient
\begin{equation*}
S_\CT(u_\CT,v_\CT) = \sum_{M\in \CM} \vartheta_M \int_M \ka_\CM \left( \nabla u_\CT \right)  \ka_\CM \left( \nabla v_\CT \right)
\end{equation*}
with
\begin{equation}\label{E:cs-lps2}
\vartheta_M  \le c_S \normL[\infty]{M}{\Va} \hdiamK[M]\quad \text{for all }M \in \CM.
\end{equation}
Here, we used the notation $I - \ka_\CM$ for  the $L^2$-projection onto an appropriate discontinuous projection space $D(\CM)$ living on the partition $\CM$ and $\bar{\Va}_{\CM}$ for a piecewise constant approximation of $\Va$ on $\CM$. The formulas for the upper bounds of $\vartheta_M$ have been discussed in detail in \cite{Kno09}. In \cite{MST07} it was shown that a local inf-sup condition
between ansatz and projection space plays an essential role in the error analysis. In the following we give some examples, for which this inf-sup condition is satisfied. 

We start with the two-level approach for which the partition $\CT$ into $d$-simplices is generated from the partition  $\CM$ into $d$-simplices by connecting the barycenter of each $M\in\CM$ with its vertices. Then, the pairs
$V(\CT) = P_D^{r,0}(\CT)$, $D(\CM) = P^{r-1,-1}(\CM)$ satisfy the inf-sup condition needed  
\cite[Lemma~3.1]{MST07}. Now let the partition $\CT$ into parallelepipeds be generated from
the partition $\CM$ into parallelepipeds by subdividing the corresponding reference cube
into $2^d$ congruent subcubes. Then, the pairs 
$V(\CT) = Q_D^{r,0}(\CT)$, $D(\CM) = Q^{r-1,-1}(\CM)$ satisfy the inf-sup condition
\cite[Lemma~3.2]{MST07}. Consequently, we could also use the pairs
$V(\CT) = Q_D^{r,0}(\CT)$, $D(\CM) = P^{r-1,-1}(\CM)$ with a smaller projection space.  

Next, for the one-level approach in which $\CM=\CT$ we introduce  scaled bubble functions $b_K\in P_{d+1}(K)\cap H_0^1(K)$ for a partition $\CT$ into $d$-simplices $K\in\CT$ and  $b_K\in Q_2(K)\cap H_0^1(K)$ for a partition $\CT$ into parallelepipeds $K\in\CT$. For $S\in\{P,Q\}$ we define the spaces
\begin{align*}
S_B^k(\CT)=\left\{ \vfi : \Om \to \rz : \vfi \rN_K=b_K\psi,\; \psi \in S_k(K) \;\text{ for all } K \in \CT \right\}.
\end{align*}
Then, the pairs $V(\CT) = P_D^{r,0}(\CT)+P_B^{r-1}(\CT)$, $D(\CM) = P^{r-1,-1}(\CT)$ on simplicial partitions $\CT$ and $V(\CT) = Q_D^{r,0}(\CT)+Q_B^{r-1}(\CT)$, $D(\CM) = P^{r-1,-1}(\CT))$ on parallelepipedal partitions $\CT$, respectively, satisfy the inf-sup condition \cite[Lemma~4.1 and 4.6]{MST07}. Note that on parallelepipedal partitions also the pairs   
$V(\CT) = Q_D^{r,0}(\CT)+P_B^{r-1}(\CT)$, $D(\CM) = P^{r-1,-1}(\CT)$ with the less enriched ansatz space satisfy the inf-sup condition \cite[Lemma~4.2]{MST07} needed for the local projection stabilization. 

\subsubsection{Subgrid scale approach}\label{SGS}
This approach, also called subgrid viscosity method, has been introduced in
\cite{Gue99,Gue01}, for an overview see eg. \cite{ErnGue04,RST08}. It is based on a scale separation
of the solution space $V(\CT)$ into a space of resolvable scales $X(\CT)$ and 
a space of unresolvable scales $Y(\CT)$ with  $V(\CT)=X(\CT)\oplus Y(\CT)$. Associated with the scale separation is a projection operator $\Pi_\CT:V(\CT)\rightarrow Y(\CT)$ with 
$X(\CT)=\ker (\Pi_\CT)$. 
As in the local projection scheme there are two types of stabilizing terms
\begin{align*}
S_\CT(u_\CT,v_\CT) &= \sum_{K\in \CT}\vartheta_K \int_K \left( \cdt \Pi_\CT(u_\CT) \right)
 \left( \cdt \Pi_\CT(v_\CT) \right),\\
 S_\CT(u_\CT,v_\CT) &= \sum_{K\in \CT} \vartheta_K \int_K  \nabla\Pi_\CT  \left(u_\CT \right) 
\nabla\Pi_\CT  \left(v_\CT \right)
 \end{align*}
with the corresponding conditions \eqref{E:cs-lps1} and \eqref{E:cs-lps2}, resp.\ for the stabilization parameters.

A typical example for spaces of resolvable and unresolvable scales on triangular partitions 
are $X(\CT)=P_D^{r,0}(\CT)$ and $Y(\CT)=P_B^{r-1}(\CT)$, $r=1,2$. One can design also subgrid scale schemes in a two-level context by setting $X(\CM)=P_D^{r,0}(\CM)$ and $V(\CT)=X(\CM)\oplus
Y(\CT)=P_D^{r,0}(\CT)$, $r=1,2,$ cf. \cite[Section~5.5]{ErnGue04}.

\subsubsection{Continuous interior penalty method} The idea of using a penalizing  term
of the form 
\begin{equation*}
S_\CT(u_\CT,v_\CT) = \sum_{E \in \CE_\Om} \vartheta_E \int_E \jump{\cd u_\CT}  \jump{\cd v_\CT}
\end{equation*}
with
\begin{equation}\label{E:cs-cip}
\vartheta_E \le c_S \hdiamE^2 \quad\text{for all }E \in \CE_\Om
\end{equation}
goes back to \cite{DD76} and has been extended to different type of problems
in  \cite{AEB07,BraBurJohLub,Bur05,Bur07,BE07,BH04,DEV13,EG13,YanZho09}.

\subsection{A Posteriori Error Estimates}\label{P:stat-err-est}
Denote by $u \in H^1_D(\Om)$ and $u_\CT \in V(\CT)$ the unique solutions of problems \eqref{E:stat-var-pb} and \eqref{E:discr-stat-pb}, respectively. Then the error $u - u_\CT$ solves the variational problem \eqref{E:stat-var-pb} with $\ell$ replaced by the residual $R$ which, for every $v \in H^1_D(\Om)$, is defined by
\begin{equation*}
\paar{R}{v} = \paar{\ell}{v} - B(u_\CT,v).
\end{equation*}
Hence, thanks to \eqref{E:robustness}, we have the following equivalence of error and residual.

\begin{LEM}[Equivalence of error and residual]\label{L:equiv-err-res}
The primal norm of the error and the dual norm of the residual are equivalent, ie.
\begin{equation*}
c_\sharp \normD{R} \le \normE{u-u_\CT} + \normD{\cd (u - u_\CT)} \le c^\sharp \normD{R}
\end{equation*}
uniformly with respect to $\ep$ and $\be$.
\end{LEM}

Integration by parts element-wise shows that the residual admits the $L^2$-re\-pre\-sen\-ta\-tion
\begin{equation*}
\paar{R}{v} = \int_\Om r v + \int_\Si j v
\end{equation*}
with
\begin{equation*}
\begin{split}
r\bigr\rvert_K &= f + \ep \De u_\CT - \cd u_\CT - b u_\CT \quad\text{for all } K \in \CT, \\
j\bigr\rvert _E &= \begin{cases}
- \jump{\ep \Vn_E \cdot \nabla u_\CT} &\text{if }E\text{ is an interior face}, \\
g - \ep \Vn_E \cdot \nabla u_\CT &\text{if }E\text{ is a face on }\Ga_N, \\
0 &\text{if }E\text{ is a face on }\Ga_D.
\end{cases}
\end{split}
\end{equation*}

Due to the stabilization term in the discrete problem \eqref{E:discr-stat-pb}, the residual $R$ does not satisfy the Galerkin orthogonality $S^{1,0}_D(\CT) \subset \ker R$ unless $S_\CT = 0$. Instead, we have for all $v_\CT \in S^{1,0}_D(\CT)$
\begin{equation*}
\paar{R}{v_\CT} = S_\CT(u_\CT,v_\CT).
\end{equation*}
To control this consistency error, denote by $I_{\CM} : H^1_D(\Om) \to S^{1,0}_D(\CM)$ any quasi-interpolation operator (cf. eg. \cite[(3.22)]{Ver13}) which satisfies, for all elements $M \in \CM$ and all faces $F$ thereof, the local error estimates (cf. eg. \cite[Prop. 3.33]{Ver13})
\begin{equation}\label{E:intpol-est}
\begin{aligned}
\normL{M}{v - I_{\CM} v} &\le c_1 \normL{\widetilde \om_{M}}{v} &\normL{M}{v - I_{\CM} v} &\le c_2 \hdiamK[M] \normL{\widetilde \om_{M}}{\nabla v} \\
\normL{M}{\nabla (v - I_{\CM} v)} &\le c_3 \normL{\widetilde \om_{M}}{\nabla v} &\normL{F}{v - I_{\CM} v} &\le c_4 \hdiamE[F]^{\halb} \normL{\widetilde \om_{F}}{\nabla v}
\end{aligned}
\end{equation}
where $\widetilde \om_{M}$ and $\widetilde \om_{F}$ denote the union of all elements in $\CM$ sharing at least a point with $M$ and $F$, respectively. The adjoint operator of $I_{\CM}$ is denoted by $I_{\CM}^\ast$ and is, for all $\vfi \in H^{-1}(\Om)$ and all $v \in H^1_D(\Om)$, defined by
\begin{equation*}
\paar{I_{\CM}^\ast \vfi}{v} = \paar{\vfi}{I_{\CM} v}.
\end{equation*}

With this notation, the above properties and \cite[Thm. 3.57, 3.59 and \S 3.8.4]{Ver13} yield the following upper and lower bounds for the dual norm of the residual.

\begin{LEM}[Bounds for the residual]\label{L:est-residual}
Define for every element $K \in \CT$ the residual a posteriori error indicator $\eta_K$ by
\begin{equation}\label{E:indicator}
\begin{split}
\eta_K &= \left\{ \mdiam{K}^2 \normL{K}{f_\CT + \ep \Delta u_\CT - \Va_\CT \cdot \nabla u_\CT - b_\CT u_\CT}^2 \phantom{\sum_{E \in \CE_K \cap \CE_{\Ga_N}}} \right. \\
&\quad\quad + \halb \sum_{E \in \CE_{K,\Om}} \ep^{- \halb} \mdiam{E} \normL{E}{\jump{\ep \Vn_E \cdot \nabla u_\CT}}^2 \\
&\quad\quad \left. + \sum_{E \in \CE_{K,\Ga_N}} \ep^{- \halb} \mdiam{E} \normL{E}{g_\CE - \ep \Vn_E \cdot \nabla u_\CT}^2 \right\}^\halb
\end{split}
\end{equation}
and the data error indicator $\theta_K$ by
\begin{equation}\label{E:data-error}
\begin{split}
\theta_K &= \left\{ \mdiam{K}^2 \normL{K}{f - f_\CT + (\Va_\CT - \Va) \cdot \nabla u_\CT + (b_\CT - b) u_\CT}^2 \phantom{\sum_{E \in \CE_K \cap \CE_{\Ga_N}}} \right. \\
&\quad\quad \left. + \sum_{E \in \CE_{K,\Ga_N}} \ep^{- \halb} \mdiam{E} \normL{E}{g - g_\CE}^2 \right\}^\halb,
\end{split}
\end{equation}
where
\begin{equation*}
\mdiam{\om} = \begin{cases}
\min \left\{ \ep^{- \halb} \diam(\om), \be^{- \halb} \right\} &\text{if } \be > 0, \\
\ep^{- \halb} \diam(\om) &\text{if } \be = 0
\end{cases}
\end{equation*}
and where $f_\CT \in S^{k,-1}(\CT)$, $\Va_\CT \in S^{k,-1}(\CT)^d$, $b_\CT \in S^{k,-1}(\CT)$, and $g_\CE \in S^{k,-1}(\CE)$ are approximations of the data $f$, $\Va$, $b$, and $g$, respectively. Then the dual norm of the residual can be bounded from above by
\begin{equation*}
\normD{R} \le c^\ast
\left\{\sum_{K \in \CT} \left[\eta_K^2 + \theta_K^2\right] \right\}^{\halb} + c^\sharp \normD{I_{\CM}^\ast R}
\end{equation*}
and from below by
\begin{equation*}
\left\{\sum_{K \in \CT} \eta_K^2\right\}^{\halb} \le c_\ast
\left[\normD{R} + \left\{\sum_{K \in \CT} \theta_K^2\right\}^{\halb} \right].
\end{equation*}
 All constants are independent of $\ep$ and $\be$; the constant $c^\sharp$ only depends on the quantity $c_b$, the constants $c^\ast$ and $c_\ast$ depend on the shape-parameter $\shape$ \eqref{E:shape-param} of $\CT$, the constant $c_\ast$ in addition depends on the polynomial degrees $k$ and $l$.
\end{LEM}

The above estimates are not yet practical since they still contain the consistency error
\begin{equation*}
\normD{I_{\CM}^\ast R} = \sup_{v \in H^1_D(\Om) \setminus \{ 0 \}} \frac{\paar{R}{I_{\CM} v}}{\normE{v}} = \sup_{v \in H^1_D(\Om) \setminus \{ 0 \}} \frac{S_\CT(u_\CT,I_{\CM} v)}{\normE{v}}.
\end{equation*}
In the next lemmas we will bound this quantity in terms of the error indicator and data errors for the stabilization schemes of the previous section.

\begin{LEM}[Consistency error of the streamline diffusion method]\label{L:sdfem}
The consistency error of the streamline diffusion method is bounded by
\begin{equation*}
\normD{I_{\CT}^\ast R} \le c \left\{ \sum_{K \in \CT} \left[ \eta_K^2 + \theta_K^2 \right] \right\}^\halb.
\end{equation*}
The constant $c$ only depends on the constants $c_S$ and $c_1, \ldots, c_4$ of equations \eqref{E:cs-sdfem} and \eqref{E:intpol-est} with $\CM=\CT$, the constant $c_I$ of equation \eqref{E:inv-est} below, and on the shape-parameter $\shape[\CT]$ \eqref{E:shape-param} of $\CT$.
\end{LEM}

\begin{BEW}
For every $v \in H^1_D(\Om)$ and every $K \in \CT$, the inverse estimate
\begin{equation}\label{E:inv-est}
\normL{K}{\nabla w_{\CT}} \le c_I \hdiamK^{-1} \normL{\widetilde \om_{K}}{w_{\CT}} \quad\forall w_{\CT} \in S^{1,0}_D(\CT)
\end{equation}
and the interpolation error estimates \eqref{E:intpol-est} imply
\begin{equation}\label{E:est-conv-der}
\begin{split}
\normL{K}{\cd (I_{\CT} v)} &\le \normL[\infty]{K}{\Va} \min \{ (1+c_3) \ep^{- \halb} \,,\, c_I (1+c_1) \hdiamK^{-1}\be^{- \halb} \}  \normE[\widetilde \om_{K}]{v}\\
&\le \max \{ 1+c_3 \,,\, c_I(1+ c_1) \} \hdiamK^{- 1} \normL[\infty]{K}{\Va} \mdiam{K} \normE[\widetilde \om_{K}]{v}.
\end{split}
\end{equation}
This estimate, assumption \eqref{E:cs-sdfem}, and the Cauchy-Schwarz inequality for integrals and sums prove the lemma.
\end{BEW}

\begin{LEM}[Consistency error of the local projection scheme]\label{L:lps}
Consider one of the variants of the local projection stabilization described in Subsection~\ref{LPS}.
In case of a partition $\CM$ of $\Omega$ into parallelepipeds, we assume that the order $r$ of the ansatz space satisfies $r\ge d$.
Then the consistency error of the local projection scheme vanishes
\begin{equation*}
\normD{I_{\CM}^\ast R} =0.
\end{equation*}
\end{LEM}
 
\begin{BEW} First we consider partitions $\CM$ into $d$-simplices (one- and two-level approach).
Since $I_{\CM} v \in P_D^{1,0}(\CM)$, we have $ \nabla (I_{\CM} v) \in P^{0,-1}(\CM)$.
Yet, $P^{0,-1}(\CM)$ is a subspace of the projection space $P^{r-1,-1}(\CM)$ which entails
\begin{equation}\label{ker}
\begin{split}
\ka_\CM(\nabla(I_\CM v))&=0,\\
\ka_\CM\left( \cdm (I_\CM v) \right)&=\bar{\Va}_{\CM}\cdot\ka_\CM 
(\nabla (I_{\CM} v))  = 0.
 \end{split}
\end{equation}
For partitions $\CM$ into parallelepipeds, we have  $I_{\CM} v \in Q_D^{1,0}(\CM)$ and 
$ \nabla (I_{\CM} v) \in P^{d-1,-1}(\CM)$. Since
$P^{d-1,-1}(\CM)\subset P^{r-1,-1}(\CM) \subset Q^{r-1,-1}(\CM)$ holds true for $r\ge d$,  
$ \nabla (I_{\CM} v)$ belongs to the kernel of $\ka_\CM$ and we again obtain \eqref{ker}.
\end{BEW}

\begin{LEM}[Consistency error of the subgrid scale approach]\label{L:sgs}
Assume that in the subgrid scale approaches described in Subsection~\ref{SGS} 
the space of resolvable scales satisfies $P_D^{1,0}(\CT)\subset X(\CT)$ or
$Q_D^{1,0}(\CT)\subset X(\CT)$ in the one-level version and $P_D^{1,0}(\CM)\subset X(\CM)$ or $Q_D^{1,0}(\CM)\subset X(\CM)$ in the two-level version, respectively. Then the consistency error vanishes
\begin{equation*}
\normD{I_{\CM}^\ast R} =0.
\end{equation*}
\end{LEM}

\begin{BEW}
Under these assumptions we have $I_\CM v\in \ker \Pi_\CM$.  
\end{BEW}

\begin{LEM}[Consistency error of the continuous interior penalty method]\label{L:cip}
Assume that $k = 1$, $V(\CT)=P_D^{1,0}$, and that the approximations $b_\CT$ and $f_\CT$ are contained in $V(\CT)$. For every element $K$ in $\CT$ set
\begin{equation}\label{E:data-error-cip}
\Theta_{\mathop{cip},K} = \mdiam{K} \normL{K}{\left( \Va - \Va_\CT \right) \cdot \nabla u_\CT} + \mdiam{K} \hdiamK \normL[\infty]{K}{\nabla \Va} \normL{K}{\nabla u_\CT}.
\end{equation}
Then the consistency error of the continuous interior penalty method is bounded by
\begin{equation*}
\normD{I_{\CT}^\ast R} \le c \left\{ \sum_{K \in \CT} \left[ \eta_K^2 + \Theta_{\mathop{cip},K}^2 \right] \right\}^\halb.
\end{equation*}
The constant $c$ only depends on the constants $c_S$, $c_1, \ldots, c_4$, $c_I$, $\widehat c_I$, and $c_{\mathop{tr}}$ of equations \eqref{E:cs-cip}, \eqref{E:intpol-est}, \eqref{E:inv-est}, \eqref{E:inv-est-face}, and \eqref{E:trace-inequ}, respectively and on the shape-parameter $\shape[\CT]$ \eqref{E:shape-param} of $\CT$.
\end{LEM}

\begin{BEW}
Since $k = 1$, we have $\De u_\CT = 0$ element-wise. Since $b_\CT$ and $f_\CT$ are supposed to be continuous, this implies for every interior face $E$ in $\CE_\Om$
\begin{equation*}
\jump{\cd u_\CT} = \jump{- \ep \De u_\CT + \Va_\CT \cdot \nabla u_\CT + b_\CT u_\CT - f_\CT} + \jump{\left( \Va - \Va_\CT \right) \cdot \nabla u_\CT}.
\end{equation*}
As usual, denote by $\om_E = K_1 \cup K_2$ the union of the two elements sharing an interior face $E$. Then, the above identity and the inverse estimate
\begin{equation}\label{E:inv-est-face}
\normL{E}{w_\CT} \le \widehat c_I \hdiamE^{- \halb} \normL{\om_E}{w_\CT} \quad\forall w_\CT \in S^{k,-1}(\CT), E \in \CE_\Om
\end{equation}
yield for every interior face $E \in \CE_\Om$
\begin{equation*}
\begin{split}
\normL{E}{\jump{\cd \left( I_\CT v \right)}} &\le \normL[\infty]{\om_E}{\Va} \normL{E}{\jump{\nabla \left( I_\CT v \right)}} \le \widehat c_I \hdiamE^{- \halb} \normL[\infty]{\om_E}{\Va} \normL{\om_E}{\nabla \left( I_\CT v \right)} \\
&\le \overline c \hdiamE^{- \frac32} \mdiam{E} \normL[\infty]{\om_E}{\Va} \normE[\om_E]{v}
\end{split}
\end{equation*}
and
\begin{equation*}
\begin{split}
\normL{E}{\jump{\cd u_\CT}} &\le \normL{E}{ \jump{- \ep \De u_\CT + \Va_\CT \cdot \nabla u_\CT + b_\CT u_\CT - f_\CT}} \\
&\quad+ \normL{E}{\jump{\left( \Va - \Va_\CT \right) \cdot \nabla u_\CT}} \\
&\le \widehat c_I \hdiamE^{- \halb} \normL{\om_E}{- \ep \De u_\CT + \Va_\CT \cdot \nabla u_\CT + b_\CT u_\CT - f_\CT} \\
&\quad+ \normL{E}{\jump{\left( \Va - \Va_\CT \right) \cdot \nabla u_\CT}}.
\end{split}
\end{equation*}
To bound the second term on the right-hand side of the last inequality, we observe that \cite[Prop. 3.5, Rem. 3.6]{Ver13} yields the trace inequality
\begin{equation}\label{E:trace-inequ}
\normL{E}{\jump{\vfi}} \le c_{\mathop{tr}} \left( \hdiamE^{- \halb} \normL{\om_E}{\vfi} + \hdiamE^{\halb} \normL{\om_E}{\nabla \vfi} \right)
\end{equation}
for every function $\vfi \in L^2(\om_E)$ with $\vfi\rvert_{K_i} \in H^1(K_i)$, $i = 1,2$. Taking again into account that $k=1$ this implies
\begin{equation*}
 \normL{E}{\jump{\left( \Va - \Va_\CT \right) \cdot \nabla u_\CT}} \le \hdiamE^{- \halb} \normL{\om_E}{\left( \Va - \Va_\CT \right) \cdot \nabla u_\CT} + \hdiamE^{\halb} \normL[\infty]{\om_E}{\nabla \Va} \normL{\om_E}{\nabla u_\CT}.
\end{equation*}
Combined with the definition of the stabilization term $S_\CT$ and of the data error $\Theta_{\mathop{cip},K}$, assumption \eqref{E:cs-cip}, and the Cauchy-Schwarz inequality for integrals and sums, this proves the bound for $\normD{I_\CT^\ast R}$.
\end{BEW}

\begin{BEM}
Replacing the stabilizing term by
\begin{equation*}
S_\CT(u_\CT,v_\CT) = \sum_{E \in \CE_\Om} \vartheta_E \int_E \jump{\Va_\CT\cdot\nabla u_\CT}  \;\jump{\Va_\CT\cdot\nabla  v_\CT}
\end{equation*}
the data error $\Theta_{cip,K}$ in Lemma~\ref{L:cip} can be omitted.
\end{BEM}
Lemmas \ref{L:equiv-err-res} -- \ref{L:cip} prove the following a posteriori error estimates.

\begin{THM}[Robust a posteriori error estimates]\label{T:error-est}
The error between the solutions $u$ and $u_\CT$ of problems \eqref{E:stat-var-pb} and \eqref{E:discr-stat-pb} can be bounded from above by
\begin{equation*}
\normE{u - u_\CT} + \normD{\cd (u - u_\CT)} \le c^\flat \left\{ \sum_{K \in \CT} \left[ \eta_K^2 + \theta_K^2 + \si_{\mathop{cip}} \Theta_{\mathop{cip},K}^2 \right] \right\}^\halb
\end{equation*}
and from below by
\begin{equation*}
\left\{ \sum_{K \in \CT} \eta_K^2 \right\}^\halb \le c_\flat \left[ \normE{u - u_\CT} + \normD{\cd (u - u_\CT)} + \left\{ \sum_{K \in \CT} \theta_K^2 \right\}^\halb \right].
\end{equation*}
Here, the assumptions for each scheme are the same as for the corresponding 
Lemmas~{\rm \ref{L:sdfem}--\ref{L:cip}}, the error indicator $\eta_K$ and the data errors $\theta_K$ and $\Theta_{\mathop{cip},K}$ are defined in equations \eqref{E:indicator}, \eqref{E:data-error}, and \eqref{E:data-error-cip}, respectively; the parameter $\si_{\mathop{cip}}$ equals~$1$ for the 
continuous interior penalty scheme and vanishes for the other discretizations.  The above error estimates are robust in the sense that the constants $c^\flat$ and $c_\flat$ are independent of the parameters $\ep$ and $\be$.
\end{THM}
%
%
\section{Non-Stationary Convection-Diffusion Equations}\label{P:non-stat-problem}
\subsection{Variational Problem}\label{P:nonstat-var-pb}
In this section, we extend the results of the previous section to the non-stationary con\-vec\-tion-dif\-fu\-sion equation
\begin{equation}\label{E:nonstat-conv-diff}
\begin{aligned}
\partial_t u - \ep \De u + \cd u + b u &= f &&\text{in } \Om \times (0,T] \\
u &= 0 &&\text{on } \Ga_D \times (0,T] \\
\ep \frac{\partial u}{\partial n} &= g &&\text{on } \Ga_N \times (0,T] \\
u(\cdot,0) &= u_0 &&\text{in } \Om
\end{aligned}
\end{equation}
in a bounded space-time cylinder with a polygonal cross-section $\Om \subset \rz^d$, $d \ge 2$, having a Lipschitz boundary $\Ga$ consisting of two disjoint parts $\Ga_D$ and $\Ga_N$. The final time $T$ is arbitrary, but kept fixed in what follows. We assume that the data satisfy the following conditions similar to assumptions (S1)--(S4):
\begin{enumerate}
\renewcommand{\labelenumi}{(T\arabic{enumi})}
\item $0 < \ep \ll 1$,
\item $f \in C(0,T;L^2(\Om))$, $g \in C(0,T;L^2(\Ga_N))$, $\Va \in C(0,T;W^{1,\infty}(\Om)^d)$, \\
$b \in C(0,T;L^\infty(\Om))$, $u_0 \in L^2(\Om)$,
\item there are two constants $\be \ge 0$ and $c_b \ge 0$, which do not depend on $\ep$, such that $-\halb \Div \Va + b \ge \be$ in $\Om \times (0,T]$ and $\normL[\infty]{\Om}{b} \le c_b \be$ in $(0,T]$,
\item the Dirichlet boundary $\Ga_D$ has positive $(d-1)$-dimensional Hausdorff measure and includes the inflow boundary $\bigcup_{0 < t \le T} \left\{x \in \Ga: \Va(x,t) \cdot \Vn(x) < 0 \right\}$.
\end{enumerate}
The first assumption of course means that we are interested in the convection-dominated regime. At the expense of more technical arguments and additional data oscillations, the second assumption can be replaced by slightly weaker conditions concerning the temporal regularity. The third assumption allows us to simultaneously handle the case of a non-vanishing reaction term and the one of absent reaction. If $b = 0$ we again set $\be = c_b = 0$.

For the variational formulation of problem \eqref{E:nonstat-conv-diff}, we denote, for every pair $a < b$ of real numbers, by $W(a,b)$ the space of all functions $v \in L^2(a,b;H^1_D(\Om))$ having their weak temporal derivative $\partial_t v$ in $L^2(a,b;H^{-1}(\Om))$ and introduce the norms
\begin{equation*}
\begin{split}
\norm[L^\infty(a,b;L^2)]{v} &= \esssup_{a < t < b} \normL{\Om}{v(\cdot,t)}, \\
\norm[L^2(a,b;H^1)]{v} &= \left\{ \int_a^b \normE{v(\cdot,t)}^2 \dx t \right\}^\halb, \\
\norm[L^2(a,b;H^{-1})]{v} &= \left\{ \int_a^b \normD{v(\cdot,t)}^2 \dx t \right\}^\halb.
\end{split}
\end{equation*}
As in \cite[\S 6.2]{Ver13}, the variational formulation of problem \eqref{E:nonstat-conv-diff} then consists in finding $u \in W(0,T)$ such that $u(\cdot,0) = u_0$ in $L^{2}(\Om)$ and such that, for almost every $t \in (0,T)$ and all $v \in H^1_D(\Om)$,
\begin{equation}\label{E:nonstat-var-pb}
\paar{\partial_t u}{v} + B(u,v) = \paar{\ell}{v}.
\end{equation}
Here the bilinear form $B$ and the linear functional $\ell$ are as in \eqref{E:def-B-ell}.

\subsection{Discretization}\label{P:nonstat-discr-pb}
For the space-time discretization of problem \eqref{E:nonstat-conv-diff}, we consider partitions $\CI = \left\{ [t_{n-1},t_n] : 1 \le n \le N_\CI \right\}$ of the time-interval $[0,T]$ into sub-intervals satisfying $0 = t_0 < \ldots < t_{N_\CI} = T$. For every $n$ with $1 \le n \le N_\CI$, we denote by $I_n = [t_{n-1},t_n]$ the $n$-th sub-interval and by $\tau_n = t_n - t_{n-1}$ its length. With every intermediate time $t_n$, $0 \le n \le N_\CI$, we associate an admissible, affine-equivalent, shape-regular partition $\CT_n$ of $\Om$ and a corresponding finite element space $V(\CT_n)$. In addition to the conditions of Section \ref{P:stat-problem}, the partitions $\CI$ and $\CT_n$ and the spaces $V(\CT_n)$ must satisfy the following assumptions.
\begin{itemize}
\item For every $n$ with $1 \le n \le N_\CI$ there is an affine-equivalent, admissible, and shape-regular partition $\widetilde\CT_n$ such that it is a refinement of both $\CT_n$ and $\CT_{n-1}$ and such that
\begin{equation*}
\shape[\widetilde \CT_n,\CT_n] = \max_{\phantom{\widetilde \CT_n}1 \le n \le N_\CI\phantom{\widetilde \CT_n}} \max_{K \in \widetilde \CT_n}
\max_{\phantom{\widetilde \CT_n}K^\prime \in \CT_n; K \subset K^\prime\phantom{\widetilde \CT_n}} \frac{\hdiamK[K^\prime]}{\hdiamK}
\end{equation*}
is of moderate size independently of $\ep$ and $\be$ (\emph{transition condition}).
\item Each $V(\CT_n)$ consists of continuous functions which are piecewise polynomials, the degrees being at least one and being bounded uniformly with respect to all partitions $\CT_n$ and $\CI$ (\emph{degree condition}).
\end{itemize}
The transition condition is due to the simultaneous presence of finite element functions defined on different grids. Usually the partition $\CT_n$ is obtained from $\CT_{n-1}$ by a combination of refinement and of coarsening. In this case the transition condition only restricts the coarsening: it should not be too abrupt nor too strong.
\\
The lower bound on the polynomial degrees is needed for the construction of suitable quasi-interpolation operators. The upper bound ensures that the constants in inverse estimates are uniformly bounded.
\\
Notice that we do not impose any shape-condition of the form $\max_n \tau_n \le c \min_n \tau_n$.

We fix a parameter $\theta \in [0,1]$ and set
\begin{equation}\label{E:ntheta}
\begin{aligned}
f^{n\theta} &= \theta f(\cdot,t_n) + (1 - \theta) f(\cdot,t_{n-1}), &g^{n\theta} &= \theta g(\cdot,t_n) + (1 - \theta) g(\cdot,t_{n-1}), \\
\Va^{n\theta} &= \theta \Va(\cdot,t_n) + (1 - \theta) \Va(\cdot,t_{n-1}), &b^{n\theta} &= \theta b(\cdot,t_n) + (1 - \theta) b(\cdot,t_{n-1})
\end{aligned}
\end{equation}
and
\begin{equation*}
\begin{split}
B^{n\theta}(u, v) &= \int_\Om \left\{ \ep \nabla u \cdot \nabla v + \Va^{n\theta} \cdot \nabla u v + b^{n\theta} u v \right\}, \\
\paar{\ell^{n\theta}}{v} &= \int_\Om f^{n\theta} v + \int_{\Ga_N} g^{n\theta} v. 
\end{split}
\end{equation*}
The finite element discretization of problem \eqref{E:nonstat-conv-diff} then consists in finding $u^n_{\CT_n} \in V(\CT_n)$, $0 \le n \le N_\CI$, such that $u^0_{\CT_0}$ is the $L^2$-projection of $u_0$ onto $V(\CT_0)$ and such that, for $n = 1, \ldots, N_\CI$ and $U^{n\theta} = \theta u^n_{\CT_n} + (1 - \theta) u^{n-1}_{\CT_{n-1}}$,
\begin{equation}\label{E:nonstat-discr-pb}
\int_\Om \frac{1}{\tau_n} (u^n_{\CT_n} - u^{n-1}_{\CT_{n-1}}) v_{\CT_n} + B^{n\theta}(U^{n\theta},v_{\CT_n}) + S^{n\theta}_{\CT_n}(U^{n\theta},v_{\CT_n}) = \paar{\ell^{n\theta}}{v_{\CT_n}}
\end{equation}
holds for all $v_{\CT_n} \in V(\CT_n)$. Here, the stabilization terms $S^{n\theta}_{\CT_n}$ are defined as in the stationary case with the following modifications:
\begin{itemize}
\item $\Va$, $b$, and $f$ are replaced by $\Va^{n\theta}$, $b^{n\theta}$, and $f^{n\theta}$, respectively,
\item $\CT$ and $\CE$ are replaced by $\CT_n$ and $\CE_n$, respectively and $\CM$ is replaced by a macro-partition $\CM_n$ subordinate to $\CT_n$,
\item for the streamline diffusion method the first factor in the integral on $K$ also contains the term $\frac{1}{\tau_n} \left( u^n_{\CT_n} - u^{n-1}_{\CT_{n-1}} \right)$.
\end{itemize}
The discrete problem \eqref{E:nonstat-discr-pb} is the familiar $\theta$-scheme and corresponds to the explicit Euler scheme, the implicit Euler scheme, and the Crank-Nicolson scheme if $\theta = 0$, $\theta = 1$, and $\theta = \halb$, respectively. We assume that, for every $n = 1, \ldots, N_\CI$, problem \eqref{E:nonstat-discr-pb} admits a unique solution $u^n_{\CT_n}$. In the literature, as for the stationary problem, this condition usually is verified by establishing the coercivity of the bilinear form $v,w \mapsto B^{n\theta}(v,w) + S^{n\theta}_{\CT_n}(v,w) - S^{n\theta}_{\CT_n}(0,w)$ with respect to a suitable mesh-dependent norm. The additional term 
$\frac{1}{\tau_n} \left( u^n_{\CT_n} - u^{n-1}_{\CT_{n-1}} \right)$ for the streamline diffusion method 
leads to a modified mass matrix and can be handled as in \cite{AMT13}.

\subsection{A Posteriori Error Estimates}\label{P:nonstat-err-est}
Denote by $u_\CI$ the function which is continuous and piecewise linear with respect to time and which equals $u^n_{\CT_n}$ at time $t_n$. With it we associate the residual $R(u_\CI) \in L^2(0,T;H^{-1})$ by setting for all $v \in H^1_D(\Om)$
\begin{equation*}
\paar{R(u_\CI)}{v} = \paar{\ell}{v} - \paar{\partial_t u_\CI}{v} - B(u_\CI,v).
\end{equation*}
Notice, that $B$ and $\ell$ are given by \eqref{E:def-B-ell} and that $\partial_t u_\CI = \frac{1}{\tau_n} \left( u^n_{\CT_n} - u^{n-1}_{\CT_{n-1}} \right)$ on $[t_{n-1},t_n]$. With this notation, we have the following equivalence of error and residual \cite[Prop. 6.14]{Ver13}.

\begin{LEM}[Equivalence of error and residual]\label{L:equi-res-error}
The error between the solutions $u$ and $u_\CI$ of problems \eqref{E:nonstat-var-pb} and \eqref{E:nonstat-discr-pb} can be bounded from below by
\begin{equation*}
\begin{split}
\norm[L^2(0,T;H^{-1})]{R(u_\CI)} &\le c^\ast \left\{ \norm[L^\infty(0,T;L^2)]{u - u_\CI}^2 + \norm[L^2(0,T;H^1)]{u - u_\CI}^2 \right. \\
&\quad\quad\quad \left. + \norm[L^2(0,T;H^{-1})]{\partial_t (u - u_\CI) + \cd (u - u_\CI)}^2 \right\}^\halb
\end{split}
\end{equation*}
and, for every $n \in \left\{ 1, \ldots, N_\CI \right\}$, from above by
\begin{equation*}
\begin{split}
&\left\{ \norm[L^\infty(0,t_n;L^2)]{u - u_\CI}^2 + \norm[L^2(0,t_n;H^1)]{u - u_\CI}^2 \right. \\
&\quad \left. + \norm[L^2(0,t_n;H^{-1})]{\partial_t (u - u_\CI) + \cd (u - u_\CI)}^2 \right\}^\halb \\
&\quad\quad\le c_\ast\left\{ \normL{\Om}{u_0 - \pi_0 u_0}^2 + \norm[L^2(0,t_n;H^{-1})]{R(u_\CI)}^2 \right\}^\halb.
\end{split}
\end{equation*}
The constants $c^\ast$ and $c_\ast$ only depend on $c_b$ and are independent of $\ep$ and $\be$.
\end{LEM}

Next, we rewrite the residual in the form
\begin{equation*}
R(u_\CI) = R_\tau (u_\CI) + R_h (u_\CI) + R_D (u_\CI).
\end{equation*}
Here, the temporal residual $R_\tau (u_\CI) \in L^2(0,T;H^{-1})$, the spatial residual $R_h (u_\CI) \in L^2(0,T;H^{-1})$, and the temporal data residual $R_D (u_\CI) \in L^2(0,T;H^{-1})$ are defined on $(t_{n-1}, t_n]$ for all $n \in \left\{ 1, \ldots, N_\CI \right\}$ and all  $v \in H^1_D(\Om)$ by
\begin{equation}\label{E:nonstat-residuals}
\begin{split}
\paar{R_\tau (u_\CI)}{v} &= B^{n\theta}(U^{n\theta} - u_\CI,v), \\
\paar{R_h (u_\CI)}{v} &= \paar{\ell^{n\theta}}{v} - \paar{\partial_t u_\CI}{v} - B^{n\theta}(U^{n\theta},v), \\
\paar{R_D (u_\CI)}{v} &= \paar{\ell}{v} - \paar{\ell^{n\theta}}{v} - B(u_\CI,v) + B^{n\theta}(u_\CI,v).
\end{split}
\end{equation}
Since $R_D(u_\CI)$ describes temporal oscillations of the known data, the task of deriving upper and lower bounds for the $L^2(t_{n-1},t_n; H^{-1})$-norms of $R(u_\CI)$ reduces to the estimation of the corresponding norms of $R_\tau(u_\CI) + R_h(u_\CI)$. The following lemma shows that this can be achieved by estimating the contributions of $R_\tau(u_\CI)$ and $R_h(u_\CI)$ separately \cite[Lemma 6.15]{Ver13}.

\begin{LEM}[Decomposition of the residual]\label{L:decomp-R}
For every $n \in \left\{ 1, \ldots, N_\CI \right\}$ we have
\begin{equation*}
\begin{split}
&\sqrt{\frac{5}{14}} \left(1 - \frac{\sqrt{3}}{2} \right) \left\{ \norm[L^2(t_{n-1},t_n; H^{-1})]{R_\tau(u_\CI)}^2 + \norm[L^2(t_{n-1},t_n; H^{-1})]{R_h(u_\CI)}^2 \right\}^\halb \\
&\quad\le \norm[L^2(t_{n-1},t_n; H^{-1})]{R_\tau(u_\CI) + R_h(u_\CI)} \\
&\quad\quad\le \norm[L^2(t_{n-1},t_n; H^{-1})]{R_\tau(u_\CI)} + \norm[L^2(t_{n-1},t_n; H^{-1})]{R_h(u_\CI)}.
\end{split}
\end{equation*}
\end{LEM}

Irrespective of the particular stabilization scheme, the temporal residual $R_\tau(u_\CI)$ equals $\left( \theta - \frac{t - t_{n-1}}{\tau_n} \right) r^n$ on the $n$-th subinterval $[t_{n-1},t_n]$, where $r^n \in H^{-1}(\Om)$ is for all $v \in H^1_D(\Om)$ defined by
\begin{equation*}
\paar{r^n}{v} = B^{n\theta}(u^n_{\CT_n} - u^{n-1}_{\CT_{n-1}},v).
\end{equation*}
An elementary calculation \cite[Lemma 6.17]{Ver13} therefore yields the following upper and lower bounds.

\begin{LEM}[Estimates for the temporal residual]\label{L:bounds-temp-res}
For every $n \in \left\{ 1, \ldots, N_\CI \right\}$, the temporal residual can be bounded from above and from below by
\begin{equation*}
\begin{split}
&c_\sharp \sqrt{\tau_n} \left\{ \normE{u^n_{\CT_n} - u^{n-1}_{\CT_{n-1}}} + \normD{\Va^{n\theta} \cdot \nabla (u^n_{\CT_n} - u^{n-1}_{\CT_{n-1}})} \right\} \\
&\quad\le \norm[L^2(t_{n-1},t_n; H^{-1})]{R_\tau(u_\CI)} \\
&\quad\quad\le c^\sharp \sqrt{\tau_n} \left\{ \normE{u^n_{\CT_n} - u^{n-1}_{\CT_{n-1}}} + \normD{\Va^{n\theta} \cdot \nabla (u^n_{\CT_n} - u^{n-1}_{\CT_{n-1}})} \right\}.
\end{split}
\end{equation*}
The constants $c_\sharp$ and $c^\sharp$ only depend on $c_b$ and are independent of $\ep$ and $\be$.
\end{LEM}

In contrast to $\normE{u^n_{\CT_n} - u^{n-1}_{\CT_{n-1}}}$ the term $\normD{\Va^{n\theta} \cdot \nabla (u^n_{\CT_n} - u^{n-1}_{\CT_{n-1}})}$ is not suited as an error indicator since it involves a dual norm. Standard approaches bound this term by inverse estimates, if need be, combined with integration by parts. This, however, leads to estimates which incorporate a factor $\ep^{- \halb}$ and which are not robust. The idea which leads to computable robust indicators is as follows \cite[Lemma 6.18]{Ver13}:
\\
Due to the definition of the dual norm, the quantities $\normD{\Va^{n\theta} \cdot \nabla (u^n_{\CT_n} - u^{n-1}_{\CT_{n-1}})}$ equal the energy-norm of the weak solutions of suitable stationary reaction-diffusion equations. These solutions are approximated by suitable finite element functions. The error of the approximations is estimated by robust error indicators for reaction-diffusion equations.

\begin{LEM}[Estimates for the convective derivative]\label{L:est-conv-der}
For every $n \in \left\{ 1, \ldots, N_\CI \right\}$ denote by $\widetilde u^{n}_{\CT_n} \in S^{1,0}_D(\widetilde \CT_n)$ the unique solution of the discrete reaction-diffusion problem
\begin{equation*}
\ep \int_\Om \nabla \widetilde u^{n}_{\CT_n} \cdot \nabla v_{\CT_n} + \be \int_\Om \widetilde u^{n}_{\CT_n}
v_{\CT_n} = \int_\Om \Va^{n\theta} \cdot \nabla (u^n_{\CT_n} - u^{n-1}_{\CT_{n-1}}) v_{\CT_n}
\end{equation*}
for all $v_{\CT_n} \in S^{1,0}_D(\widetilde \CT_n)$. Define the error indicator $\widetilde \eta^{n}_{\CT_n}$ by
\begin{equation*}
\begin{split}
\widetilde \eta^{n}_{\CT_n} &= \left\{ \sum_{K \in \widetilde \CT_n}
\mdiam{K}^2 \normL{K}{\Va^{n\theta} \cdot \nabla (u^n_{\CT_n} - u^{n-1}_{\CT_{n-1}}) + \ep \De \widetilde u^{n}_{\CT_n} - \be \widetilde u^{n}_{\CT_n}}^2 \right. \\
&\quad\quad \left. + \sum_{E \in \widetilde \CE_{n,\Om} \cup \widetilde \CE_{n,\Ga_N}} \ep^{- \halb} \mdiam{E}
\normL{E}{\jump{\Vn_E \cdot \nabla \widetilde u^{n}_{\CT_n}}}^2 \right\}^\halb.
\end{split}
\end{equation*}
Then there are two constants $c_\dagger$ and $c^\dagger$ which only depend on the shape-parameter $\shape[\widetilde \CT_n]$ \eqref{E:shape-param} of $\widetilde \CT_n$ such that the following estimates are valid
\begin{equation*}
c_\dagger \left\{ \normE{\widetilde u^{n}_{\CT_n}} + \widetilde \eta^{n}_{\CT_n} \right\} \le \normD{\Va^{n\theta} \cdot \nabla (u^n_{\CT_n} - u^{n-1}_{\CT_{n-1}})} \le c^\dagger
\left\{ \normE{\widetilde u^{n}_{\CT_n}} + \widetilde \eta^{n}_{\CT_n} \right\}.
\end{equation*}
\end{LEM}

A comparison of equations \eqref{E:nonstat-discr-pb} and \eqref{E:nonstat-residuals} reveals that, on each interval $(t_{n-1},t_n]$ separately, the spatial residual $R_h(u_\CI)$ is the residual of a stationary problem \eqref{E:def-B-ell} with suitably modified functions $\Va$, $b$, $f$, and $g$. Hence, the results of Section \ref{P:stat-problem} yield the following upper and lower bounds for the dual norm of the spatial residual.

\begin{LEM}[Estimates for the spatial residual]\label{L:bounds-spatial-res}
For every $n \in \{ 1, \ldots, N_\CI \}$ define a spatial error indicator by
\begin{equation*}
\begin{split}
\eta^n_{\CT_n} &= \left\{ \sum_{K \in \widetilde \CT_n} \mdiam{K}^2 \normL{K}{f^{n\theta}_{\CT_n} - \frac{u^n_{\CT_n} - u^{n-1}_{\CT_{n-1}}}{\tau_n} + \ep \Delta U^{n\theta} - \Va^{n\theta}_{\CT_n} \cdot \nabla U^{n\theta} - b^{n\theta}_{\CT_n} U^{n\theta}}^2 \phantom{\sum_{\CE_{\Ga_N}}} \right. \\
&\quad\quad + \halb \sum_{E \in \widetilde \CE_{n,\Om}} \ep^{- \halb} \mdiam{E} \normL{E}{\jump{\ep \Vn_E \cdot \nabla U^{n\theta}}}^2 \\
&\quad\quad \left. + \sum_{E \in \widetilde \CE_{n,\Ga_N}} \ep^{- \halb} \mdiam{E} \normL{E}{g^{n\theta}_{\CE_n} - \ep \Vn_E \cdot \nabla U^{n\theta}}^2 \right\}^\halb
\end{split}
\end{equation*}
and spatial data errors by
\begin{align*}
\theta^n_{\CT_n} &= \left\{ \sum_{K \in \CT_n} \mdiam{K}^2 \normL{K}{f^{n\theta} - f^{n\theta}_{\CT_n} + (\Va^{n\theta}_{\CT_n} - \Va^{n\theta}) \cdot \nabla U^{n\theta} + (b^{n\theta}_{\CT_n} - b^{n\theta}) U^{n\theta}}^2 \phantom{\sum_{\CE_{\Ga_N}}} \right. \\*
&\quad\quad \left. + \sum_{E \in \CE_{n,\Ga_N}} \ep^{- \halb} \mdiam{E} \normL{E}{g^{n\theta} - g^{n\theta}_{\CE_n}}^2 \right\}^\halb, \\
\Theta^n_{\mathop{cip},\CT_n} &= \left\{ \sum_{K \in \CT_n} \mdiam{K}^2 \normL{K}{\left( \Va^{n\theta} - \Va^{n\theta}_{\CT_n} \right) \cdot \nabla U^{n\theta}}^2 + \mdiam{K}^2 \hdiamK^2 \normL[\infty]{K}{\nabla \Va^{n\theta}} \normL{K}{\nabla U^{n\theta}}^2 \right\}^\halb.
\end{align*}
Here, $U^{n\theta} = \theta u^n_{\CT_n} + (1 - \theta) u^{n-1}_{\CT_{n-1}}$ is as in \eqref{E:nonstat-discr-pb}, $f^{n\theta}$, $g^{n\theta}$, $\Va^{n\theta}$, and $b^{n\theta}$ are as in \eqref{E:ntheta}, and $f^{n\theta}_{\CT_n}$, $\Va^{n\theta}_{\CT_n}$, $b^{n\theta}_{\CT_n}$, and $g^{n\theta}_{\CE_n}$ are approximations of $f^{n\theta}$, $\Va^{n\theta}$, $b^{n\theta}$, and $g^{n\theta}$ on $\CT_n$ and $\CE_n$, respectively. Then, on every interval $(t_{n-1},t_n]$, the dual norm of the spatial residual can be bounded from above by
\begin{equation*}
\normD{R_h(u_\CI)} \le c^\flat \left\{ \left( \eta^n_{\CT_n} \right)^2 + \left( \theta^n_{\CT_n} \right)^2 
+ \si_{\mathop{cip}} \left( \Theta^n_{\mathop{cip},\CT_n} \right)^2 \right\}^\halb
\end{equation*}
and from below by
\begin{equation*}
\eta^n_{\CT_n} \le c_\flat \left[ \normD{R_h(u_\CI)} + \theta^n_{\CT_n} \right].
\end{equation*}
Here, the assumptions for each stabilized scheme are the same as for the corresponding 
Lemmas~{\rm \ref{L:sdfem}--\ref{L:cip}},
the parameter $\si_{\mathop{cip}}$ equals $1$ for the continuous interior penalty scheme and vanishes for the other discretizations. The above error estimates are robust in the sense that the constants 
$c^\flat$ and $c_\flat$ are independent of the parameters $\ep$ and $\be$.
\end{LEM}

Lemmas \ref{L:equi-res-error} -- \ref{L:bounds-spatial-res} yield the following a posteriori error estimates for the non-stationary problem.

\begin{THM}[Robust a posteriori error estimates]\label{T:nonstat-error-est}
The error between the solution $u$ of problem \eqref{E:nonstat-var-pb} and the solution $u_\CI$ of problem \eqref{E:nonstat-discr-pb} is bounded from above by
\begin{equation*}
\begin{split}
&\left\{ \norm[L^\infty(0,T;L^2)]{u - u_\CI}^2 + \norm[L^2(0,T;H^1)]{u - u_\CI}^2 \right. \\
&\quad \left. + \norm[L^2(0,T;H^{-1})]{\partial_t (u - u_\CI) + \cd (u - u_\CI)}^2 \right\}^\halb \\
&\quad\quad\le c^\ast \left\{ \normL{\Om}{u_0 - \pi_0 u_0}^2 \phantom{\sum_{n = 1}^{N_\CI}} \right. \\
&\quad\quad\quad\phantom{c^\ast \Bigl\{}
+ \sum_{n = 1}^{N_\CI} \tau_n \left[ \left(\eta^n_{\CT_n}\right)^2 + \normE{u^n_{\CT_n} - u^{n-1}_{\CT_{n-1}}}^2 + \left(\widetilde \eta^{n}_{\CT_n}\right)^2 + \normE{\widetilde u^{n}_{\CT_n}}^2 \right] \\
&\quad\quad\quad\phantom{c^\ast \Bigl\{}
+ \sum_{n = 1}^{N_\CI} \tau_n \left[ \left(\theta^n_{\CT_n}\right)^2 
+ \si_{\mathop{cip}} \left( \Theta^n_{\mathop{cip},\CT_n} \right)^2 \right] \\
&\quad\quad\quad\phantom{c^\ast \Bigl\{}
+ \norm[L^2(0,T;H^{-1})]{f - f^{n\theta} - (\Va - \Va^{n\theta}) \cdot \nabla
u_\CI - (b - b^{n\theta}) u_\CI}^2 \\
&\quad\quad\quad\phantom{c^\ast \Bigl\{}
\left. + \sum_{n = 1}^{N_\CI} \norm[L^2(t_{n-1},t_n;H^{- \halb}(\Ga_N))]{g - g^{n\theta}_{\CE_n}}^2
\right\}^\halb
\end{split}
\end{equation*}
and on each interval $(t_{n-1},t_n]$, $1 \le n \le N_\CI$, from below by
\begin{equation*}
\begin{split}
&\tau_n^\halb \left\{\left(\eta^n_{\CT_n}\right)^2 + \normE{u^n_{\CT_n} - u^{n-1}_{\CT_{n-1}}}^2 + \left(\widetilde \eta^{n}_{\CT_n}\right)^2 + \normE{\widetilde u^{n}_{\CT_n}}^2 \right\}^\halb \\
&\quad\le c_\ast \left\{ \norm[L^\infty(t_{n-1},t_n;L^2)]{u - u_\CI}^2 + \norm[L^2(t_{n-1},t_n;H^1)]{u - u_\CI}^2 \right. \\
&\quad\quad\quad+ \norm[L^2(t_{n-1},t_n;H^{-1})]{\partial_t (u - u_\CI) + \cd (u - u_\CI)}^2 \\ 
&\quad\quad\quad+ \tau_n \left(\theta^n_{\CT_n}\right)^2 \\
&\quad\quad\quad+ \norm[L^2(t_{n-1},t_n;H^{-1})]{f - f^{n\theta} - (\Va - \Va^{n\theta}) \cdot \nabla
u_\CI - (b - b^{n\theta}) u_\CI}^2 \\
&\quad\quad\quad
\left. + \norm[L^2(t_{n-1},t_n;H^{- \halb}(\Ga_N))]{g - g^{n\theta}_{\CE_n}}^2
\right\}^\halb.
\end{split}
\end{equation*}
Here, the assumptions for each scheme are the same as for the corresponding 
Lemmas~{\rm \ref{L:sdfem}--\ref{L:cip}},
the functions $\widetilde u^{n}_{\CT_n}$ and the indicators $\widetilde \eta^{n}_{\CT_n}$ are defined in Lemma \ref{L:est-conv-der}, and the quantities $\eta^n_{\CT_n}$, $\theta^n_{\CT_n}$,  and $\Theta^n_{\mathop{cip},\CT_n}$ are as in Lemma \ref{L:bounds-spatial-res}. The parameter $\si_{\mathop{cip}}$ equals $1$ for the continuous interior penalty scheme and vanishes for the other discretizations. The above error estimates are robust in the sense that the constants $c^\ast$ and $c_\ast$ are independent of the final time $T$, the viscosity $\ep$ and the parameter $\be$.
\end{THM}
%
%
\providecommand{\bysame}{\leavevmode\hbox to3em{\hrulefill}\thinspace}
\providecommand{\MR}{\relax\ifhmode\unskip\space\fi MR }
\providecommand{\MRhref}[2]{%
  \href{http://www.ams.org/mathscinet-getitem?mr=#1}{#2}
}
\providecommand{\href}[2]{#2}

\end{document}